\documentclass[12pt]{article}

\usepackage{amsmath}
\usepackage{amssymb}
\usepackage{graphicx}
\usepackage{epsfig}

\usepackage[margin=2cm]{geometry}

\newtheorem{theorem}{Theorem}
\newtheorem{lemma}[theorem]{Lemma}

\newtheorem{proposition}[theorem]{Proposition}

\def \Z {{\mathbb{Z}}}

\def \proof {\noindent{\bf Proof}\quad}

\def \md#1{{\,({\rm mod}\ #1)}}

\def\endmark{\hskip 2em$\square$\par}
\def \qed {\hfill\endmark}

\def \cay {{\rm Cay}}
\def \sym {{\rm Sym}}
\def \alt {{\rm Alt}}
\def \gl {{\rm GL}}
\def \psl {{\rm PSL}}
\def \aut {{\rm Aut}}

\title{\bf
Vertex-transitive graphs that have\\ no Hamilton decomposition}

\author{
Darryn Bryant
\thanks{
School of Mathematics and Physics,
The University of Queensland,
Qld 4072, Australia.
\texttt{db@maths.uq.edu.au},
}
\and
Matthew Dean
\thanks{
School of Mathematics and Physics,
The University of Queensland,
Qld 4072, Australia.
\texttt{mdean@uq.edu.au},
}
}

\date{ }

\begin{document}
\maketitle\thispagestyle{empty}
\def\baselinestretch{1.5}\small\normalsize

\begin{abstract}
It is shown that there are infinitely many connected vertex-transitive graphs that have no Hamilton decomposition, including infinitely many Cayley graphs of valency $6$, and including Cayley graphs of arbitrarily large valency. 
\end{abstract}

\section{Introduction}

A famous question of Lov\'asz concerns the existence of Hamilton paths in vertex-transitive graphs \cite{Lov}, and no example of a connected vertex-transitive graph with no Hamilton path is known. The related question concerning the 
existence of Hamilton cycles in vertex-transitive graphs is another interesting and well-studied problem in graph theory, see the survey \cite{KutMar}. A {\em Hamiltonian} graph is a graph containing a Hamilton cycle. Thomassen (see \cite{Ber,KutMar}) has conjectured that there are only finitely many non-Hamiltonian connected vertex-transitive graphs. On the other hand, Babai \cite{Bab1,Bab2} has conjectured that there are infinitely many such graphs. To date only five are known. These are the complete graph of order $2$, the Petersen graph, the Coxeter graph, and the two 
graphs obtained from the Petersen and Coxeter graphs by replacing each vertex with a triangle. 

For a regular graph of valency at least 4, a stronger property than the existence of a Hamilton cycle is the existence of a Hamilton decomposition.
If $X$ is a $k$-valent graph, then a {\it Hamilton decomposition} of $X$ is a set of
$\lfloor\frac k2\rfloor$ pairwise edge-disjoint Hamilton cycles in $X$.
Given the small number of non-Hamiltonian connected vertex-transitive graphs, 
and the uncertainty concerning the existence of others, it is natural to ask how many connected vertex-transitive graphs have no Hamilton decomposition. 

Mader \cite{Mad} showed that a connected $k$-valent vertex-transitive graph 
is $k$-edge-connected. So for any connected vertex-transitive graph $X$, 
there is no obvious obstacle to the existence of a Hamilton decomposition of $X$.
Indeed, Wagon has conjectured that with a handful of small exceptions,
every connected vertex-transitive graph has a Hamilton decomposition, see \cite{Wei}. 
As discussed below, there is a lot of evidence to support this conjecture. 
However, in this paper we show that there are in fact infinitely many connected vertex-transitive graphs that have no Hamilton decomposition, 
including infinitely many connected $6$-valent Cayley graphs, 
and including Cayley graphs of arbitrarily large valency.

As far as we are aware, 
there are six previously known examples of connected vertex-transitive graphs 
that have no Hamilton decomposition, and none of these is a Cayley graph. 
Firstly, there are the four non-Hamiltonian $3$-valent graphs mentioned above. Secondly, Kotzig \cite{Kot} has shown that 
a $3$-valent graph has a Hamilton cycle if and only if its line graph has a Hamilton decomposition. Thus, 
the line graphs of the four known non-Hamiltonian connected $3$-valent vertex-transitive graphs
are $4$-valent graphs that have no Hamilton decomposition. However, of these, only the line graphs of the Petersen and Coxeter graphs are vertex-transitive. 
Wagon \cite{Wei} has verified that 
every other connected vertex-transitive graph of order at most $31$ has a Hamilton decomposition. We have independently verified this, 
using McKay and Royle's list of vertex-transitive graphs that is available online 
(also see \cite{McKRoy}).

Poto\u cnik, Spiga and Verret have found that there are $4,820$ connected $4$-valent 
arc-transitive graphs with at most 640 vertices \cite{PotSpiVer}, and McKay
has shown by computation that all of these have Hamilton decompositions, except the line graphs of the Petersen and Coxeter graphs. 
Alspach and Rosenfeld \cite{AlsRos} have asked whether every prism over a connected
$3$-valent Hamiltonian graph has a Hamilton decomposition. The prism over a graph $X$ is the cartesian product of $X$ and the complete graph of order $2$. 
McKay has shown by computation that the prism 
over a connected $3$-valent vertex-transitive graph of order at most $500$ has a Hamilton decomposition (prisms over vertex-transitive graphs are vertex-transitive).

Existence of Hamilton decompositions of vertex-transitive graphs has been established in many other cases. Alspach \cite{Als3} showed that every connected vertex-transitive graph of order $2p$, where $p\equiv 3\md 4$ is prime, has a Hamilton decomposition. More recently, it has been proved that all connected vertex-transitive graphs of order $p$ or $p^2$, where $p$ is prime, have Hamilton decompositions \cite{AlsBryKre}. Every such graph is in fact a connected Cayley graph on an abelian group, and a long-standing conjecture of Alspach \cite{Als1,Als2} is that every connected Cayley graph on an abelian group has a Hamilton decomposition. 
This conjecture has been verified for graphs with valency at most 5 \cite{AlsHeiLiu,BerFavMah,CheQui,FanLicLiu}, and in many cases for valency 6 \cite{Dea1,Dea2,Wes1,Wes2,WesLiuKre}.
Also, Liu \cite{Liu1,Liu2,Liu3} has proved strong results on the problem in cases where restrictions are placed on the connection set of the graph. We show that Alspach's conjecture does not extend to Cayley graphs on non-abelian groups by exhibiting several infinite families of connected Cayley graphs that have no Hamilton decomposition.

Hamilton decompositions of general graphs, not necessarily vertex-transitive, have been studied extensively, see the survey \cite{Gou}. A very well-known conjecture on Hamilton decompositions is due to Nash-Williams \cite{Nas}. The slightly strengthened version of his conjecture, due to Jackson \cite{Jac}, states that every connected $k$-valent graph of order at most $2k+1$ has a Hamilton decomposition. This conjecture has recently been proved for all sufficiently large $k$ by
Csaba, K\"uhn, Lo, Osthus and Treglown \cite{CsaKuhLoOstTre}.
Another result, due to Gr\"unbaum and Malkevitch \cite{GruMal}, is that there exist $4$-valent $4$-connected graphs that have no Hamilton decomposition, and moreover that there exist planar graphs with this property. We make use of one of the main ideas from their paper. There are also two papers by Pike \cite{Pik1,Pik2} that concern Hamilton decompositions, and in particular contain some questions on the existence of Hamilton decompositions of vertex-transitive graphs. 

\section{Preliminaries}

For $3$-valent arc-transitive graphs, we use a notation that is consistent with common usage, such as in \cite{ConDob}. 
A $3$-valent arc-transitive graph is denoted by $F$, followed by its order, 
followed by a letter ($A$, $B$, $C$, and so on) when there is more than one $3$-valent arc-transitive graph 
of a given order. For example, the Petersen graph is denoted $F10$, and the two $3$-valent arc-transitive
graphs of order $20$, the Dodecahedron graph and the Desargues graph, are denoted by $F20A$ and $F20B$. 
The common names of the twelve connected $3$-valent arc-transitive graphs of orders at least $8$ and at most $32$ are given in the following table. 

\begin{table}[ht]
\centering 
\begin{tabular}{|ll|ll|ll|} 
\hline 
$F8$: & 3-cube graph & $F10$: & Petersen graph & $F14$: & Heawood graph \\ 
\hline
$F16$: & M\"obius-Kantor graph & $F18$: & Pappus graph & $F20A$: & Dodecahedron graph\\ 
\hline
$F20B$: & Desargues graph & $F24$: & Nauru graph & $F26$: & F26A graph\\ 
\hline
$F28$: & Coxeter graph & $F30$: & Tutte-Coxeter graph & $F32$: & Dyck graph\\ 
\hline 
\end{tabular}
\end{table}

We will be dealing with multigraphs, and we need to take some care with the notation used. 
Any graph is understood to be simple, and we use the term multigraph whenever there are distinct edges with the same endpoints. None of our graphs or multigraphs have loops. 
In any graph, we use $\{x,y\}$ to denote the unique edge with endpoints $x$ and $y$. 
Similarly, we use $(x,y)$ to denote the unique arc from $x$ to $y$. 
For any given graph $X$, the multigraph denoted by $mX$ has the same vertices as $X$, 
and has $m$ distinct edges $\{x,y\}_0,\{x,y\}_1,\ldots,\{x,y\}_{m-1}$ joining $x$ and $y$ 
for each edge $\{x,y\}$ in $X$. In $mX$, we distinguish $m$ arcs $(x,y)_0,(x,y)_1,\ldots,(x,y)_{m-1}$ 
for each arc $(x,y)$ in $X$, and associate the two arcs $(x,y)_i$ and $(y,x)_i$ of $mX$ with the edge 
$\{x,y\}_i$ of $mX$. 
 
Let $X$ be a non-empty regular graph of valency $k$ and order $n$, and let $m$ be a positive integer.
We define $K(mX)$ as follows. 
The vertices of $K(mX)$ are the arcs of $mX$. 
For each vertex $v$ of $mX$, there is a complete subgraph of $K(mX)$ on the $km$
arcs emanating from $v$. We refer to this complete subgraph of $K(mX)$ as the complete subgraph 
associated with $v$. 
Also, for each edge $\{x,y\}_i$ in $mX$, there is an edge in $K(mX)$ 
joining $(x,y)_i$ and $(y,x)_i$, and we associate the edge $\{(x,y)_i,(y,x)_i\}$ of $K(mX)$ with 
the edge $\{x,y\}_i$ of $mX$.
This is all the edges of $K(mX)$. When $m=1$ we may write just $K(X)$ rather than $K(1X)$.
 
It should be apparent that $K(mX)$ is isomorphic to the graph obtained from $mX$ by replacing 
each vertex of $mX$ with a complete graph of order $km$. 
Observe that $K(mX)$ is a regular graph of valency $km$ and order $kmn$, 
and that $K(X)$ is connected if and only if $X$ is connected. 
In \cite{AlsDob}, various properties of these graphs are proved (for the case $m=1$).
In particular, it is shown that if $m=1$ and $X$ is a connected vertex-transitive graph 
of valency $k\geq 3$, then $K(mX)$ is vertex-transitive if and only if 
$X$ is arc-transitive. Using our above definition of $K(mX)$, 
it is easy to see that this result is in fact true for all $k\geq 1$ and for all 
$m\geq 1$.

\begin{lemma}\label{GKVTiffGAT}
If $X$ is a non-empty regular graph and $m$ is a positive integer, 
then the graph $K(mX)$ is vertex-transitive if and only if $X$ is arc-transitive. 
\end{lemma}

\begin{lemma}\label{GHDiffGKHD}
Let $X$ be a regular graph and let $m$ be a positive integer.
The graph $K(mX)$ has a Hamilton decomposition if and only if $mX$ has a Hamilton decomposition.
\end{lemma}

\proof
Let $k$ be the valency of $X$ and let $t=\lfloor\frac{km}2\rfloor$.
For each vertex $v$ in $mX$,
let $X_v$ be the complete subgraph of $K(mX)$ associated with the vertex $v$,
and let $E_v$ be the set of 
edges of $K(mX)$ having exactly one endpoint in $X_v$. Equivalently, $E_v$ is the set of edges of 
$K(mX)$ associated with the edges of $mX$ that are incident on $v$. 
Observe that $|E_v|=km$. 

First suppose $K(mX)$ has a Hamilton decomposition $\{Y_1,Y_2,\ldots,Y_t\}$. 
For $v\in V$ and $1\leq j\leq t$, the number of edges of $Y_j$ in $E_v$ is positive 
and even. Since $|E_v|=km$, it follows that this number is $2$. 
Hence, if we contract the edges of each $X_v$, then 
each $Y_j$ contracts to a Hamilton cycle $C_j$ in $mX$, 
and $\{C_j:1\leq j\leq t\}$ is a Hamilton decomposition of $mX$.

Now, conversely, suppose that $mX$ has a Hamilton decomposition $\{C_1,C_2,\ldots,C_t\}$.
For $1\leq j\leq t$, let 
$Z_j$ consist of the edges of $K(mX)$ that have an associated edge in $C_j$.
It is clear that each $Z_j$ can be extended to the edge set of a Hamilton 
cycle in $K(mX)$ by adding the edges of a Hamilton path in each $X_v$, such that 
each Hamilton path has the required endpoints.

If $km$ is even, then the complete graph of order $km$ 
can be decomposed into $\frac{km}2$ pairwise edge-disjoint Hamilton paths,
and in any such decomposition each vertex is an endpoint of exactly one of the Hamilton paths.
Also, if $km$ is odd, then the complete graph of order $km$ can be decomposed into $\frac{km-1}2$ 
pairwise edge-disjoint 
Hamilton paths and a matching of order $mk-1$. 
In any such decomposition each vertex of the 
matching is an endpoint of exactly one of the Hamilton paths. 
Thus, both when $km$ is even and when $km$ is odd, 
$Z_1,Z_2,\ldots,Z_t$ can be extended to a Hamilton decomposition of 
$K(mX)$.
\qed

\vspace{0.3cm}

We are interested in connected vertex-transitive graphs that have no Hamilton decomposition,
and the following immediate consequence of 
Lemmas \ref{GKVTiffGAT} and \ref{GHDiffGKHD} gives us a method for constructing them.

\begin{lemma}\label{mainconstruction}
If $X$ is a connected arc-transitive graph and $mX$ has no Hamilton decomposition, 
then $K(mX)$ is a
connected vertex-transitive graph that has no Hamilton decomposition. 
\end{lemma}

The line graph of a graph $X$ is denoted by $L(X)$.
Since $L(F10)$ and $L(F28)$ are arc-transitive and have no Hamilton decomposition, Lemma \ref{mainconstruction} tells us that $K(L(F10))$ and $K(L(F28))$ are 
$4$-valent vertex-transitive graphs that have no Hamilton decomposition. 
These two graphs are in fact Cayley graphs, and represent the first examples of connected Cayley graphs 
that are known to have no Hamilton decomposition.
The graph $K(L(F10))$ is a Cayley graph on the alternating group $\alt(5)$, and $K(L(F28))$ is a Cayley graph on the projective special linear group $\psl(2,7)$. 
This can be seen by noting 
the correspondence between the $2$-arcs of a graph $X$ and the vertices of $K(L(X))$,
that $\alt(5)$ has a regular action on the $2$-arcs of $F10$,
and that $\psl(2,7)$ has a regular action on the $2$-arcs of $F28$,
see \cite{ConNed}.

\begin{proposition}
The graphs $K(L(F10))$ and $K(L(F28))$
are connected $4$-valent Cayley graphs that have no Hamilton decomposition,
where $F10$ is the Petersen graph and $F28$ is the Coxeter graph.
\end{proposition}

\section{6-valent vertex-transitive graphs}\label{6valentSection}

For each $3$-valent arc-transitive graph $X$ of order at most $50$, we have verified by computer whether 
$2X$ has a Hamilton decomposition. If $2X$ has no Hamilton decomposition, 
then by Lemma \ref{mainconstruction}, $K(2X)$ is a $6$-valent vertex-transitive graph with 
no Hamilton decomposition. The results of our computer search give us the following proposition.
  
\begin{proposition}\label{afewmore}
The graphs $K(2F8)$, $K(2F10)$, $K(2F16)$, $K(2F18)$, $K(2F20B)$, $K(2F24)$, $K(2F28)$, $K(2F30)$,
$K(2F32)$, $K(2F40)$, $K(2F48)$ and $K(2F50)$ are connected $6$-valent vertex-transitive graphs that have no Hamilton decomposition.
\end{proposition}

We now proceed to show the existence of infinitely many connected $6$-valent vertex-transitive graphs that have no Hamilton decomposition. 
The following lemma shows that if $X$ is $3$-valent, then the existence of a Hamilton decomposition of $2X$ 
is equivalent to the existence of a perfect $1$-factorisation of $X$. 
A {\em perfect $1$-factorisation} of a $k$-valent graph is a set of $k$ pairwise edge-disjoint $1$-factors (perfect matchings) such that the union of any two of these $1$-factors is a Hamilton cycle.

\begin{lemma}\label{PerfactiffHD}
If $X$ is a $3$-valent graph, then $X$ has a perfect $1$-factorisation if and only if $2X$ has a Hamilton decomposition. 
\end{lemma}

\proof
If $\{X_1,X_2,X_3\}$ is a perfect $1$-factorisation of $X$, then 
$\{X_1\cup X_2,X_1\cup X_3,X_2\cup X_3\}$ yields a Hamilton decomposition of $2X$. 
Conversely, if $\{Y_1,Y_2,Y_3\}$ is a Hamilton decomposition of $2X$, 
and we let $X_1$ contain those edges of $X$ where 
the corresponding two edges of $2X$ are in $Y_1$ and $Y_2$, let
$X_2$ contain those edges of $X$ where the corresponding two edges of $2X$ are in $Y_1$ and $Y_3$, and let
$X_3$ contain those edges of $X$ where the corresponding two edges of $2X$ are in $Y_2$ and $Y_3$, then 
$\{X_1,X_2,X_3\}$ is a perfect $1$-factorisation of $X$.
\qed

\vspace{0.3cm}

An immediate corollary of Lemma \ref{PerfactiffHD} (combined with Lemma \ref{mainconstruction})
is that if $X$ is a $3$-valent arc-transitive graph that has 
no perfect $1$-factorisation, then $K(2X)$ is a $6$-valent vertex-transitive graph that has no Hamilton decomposition. The following result, which Laufer \cite{Lau} attributes to Kotzig \cite{Kot}, is thus 
important for us. 

\begin{theorem}\rm{(Kotzig, \cite{Kot})}\label{KotzigTheorem}
If $X$ is a regular bipartite graph of order congruent to $0\md 4$ and valency at least $3$, 
then $X$ has no perfect $1$-factorisation.
\end{theorem}

Theorem \ref{KotzigTheorem} combined with Lemmas \ref{mainconstruction} and \ref{PerfactiffHD} gives us the following theorem. 

\begin{theorem}\label{0mod4bipartite}
If $X$ is a connected bipartite $3$-valent arc-transitive graph of order congruent to 
$0\md 4$, then $K(2X)$ is a connected $6$-valent vertex-transitive graph that has no Hamilton decomposition.
\end{theorem}

All except five of the graphs in Proposition \ref{afewmore} are of the form $K(2X)$ where $X$ is a 
bipartite graph of order congruent to $0\md 4$. 
The exceptions are $K(2F10)$, $K(2F18)$, $K(2F28)$, $K(2F30)$ and $K(2F50)$.
Since 
it is known that there are infinitely many connected bipartite $3$-valent arc-transitive graphs of order congruent to 
$0\md 4$, see \cite{ConNed} for example, we have the following corollary to 
Theorem \ref{0mod4bipartite}.

\begin{theorem}\label{InfinitelyManySixRegular}
There are infinitely many $6$-valent connected vertex-transitive graphs that have no Hamilton decomposition.
\end{theorem}

Many of the graphs given by Theorem \ref{InfinitelyManySixRegular} 
are Cayley graphs. 
To see this, consider the action of $\aut(X)\times\Z_m$ on the vertices of $K(mX)$ given by 
$$(x,y)_i(g,a)=(xg,yg)_{i+a}$$
for each $(g,a)\in\aut(X)\times\Z_m$ and each vertex $(x,y)_i$ of $K(mX)$.
Here, the subscript $i+a$ is calculated in $\Z_m$. It is easily seen that $\aut(X)\times\Z_m$ is a subgroup
of $\aut(K(mX))$. Moreover, noting that the arcs of $mX$ are the vertices of $K(mX)$, we see 
that if $G$ is a subgroup of $\aut(X)$ with a regular action on the arcs of $X$, 
then $G\times \Z_m$ is a subgroup of $\aut(K(mX))$ with a regular action on the vertices of $K(mX)$. Thus, $K(mX)$ is a Cayley graph. 

In \cite{ConNed}, connected $3$-valent arc-transitive graphs which admit a regular group action on their arcs are referred to as having a {\em Type 1 action},
and it is shown that there are infinitely many such graphs that are 
bipartite and have order congruent to $0\md 4$.
Combining this with Theorem \ref{InfinitelyManySixRegular} and the discussion of the preceding paragraph we have the following result. 

\begin{theorem}\label{Cayley6regular}
There are infinitely many connected $6$-valent Cayley graphs that have no Hamilton decomposition.
\end{theorem}

We note that not all the $6$-valent connected vertex-transitive graphs with no Hamilton decomposition that we have constructed are Cayley graphs. For example, consider the graph 
$K(2F30)$. It is known that $\aut(F30)\cong\sym(6)\times\Z_2$, and so it is easily seen that 
$\aut(K(2F30))\cong\sym(6)\times\Z_2\times\Z_2^{45}$. Thus, if $K(2F30)$ is a Cayley graph, then $\sym(6)\times\Z_2\times\Z_2^{45}$ has a subgroup of order $180$
(the order of $K(2F30)$). Since $\sym(6)$ has no subgroup of order $45$, $90$ or $180$, this is not the case. It follows that $K(2F30)$ is not a Cayley graph.

\section{Vertex-transitive graphs of arbitrarily large valency}

For each positive integer $m$, the multigraph
$mF10$ is arc-transitive and has no Hamilton decomposition 
(because $F10$ is arc-transitive and non-Hamiltonian).
It thus follows from Lemma \ref{mainconstruction} that 
$K(mF10)$ is a connected $3m$-valent vertex-transitive graph of 
order $30m$ that has no Hamilton decomposition.
Similarly, $K(mF28)$ is a connected $3m$-valent vertex-transitive graph of 
order $84m$ that has no Hamilton decomposition.
Thus, there exist connected vertex-transitive graphs of 
arbitrarily large valency that have no Hamilton decomposition.
 
We now give two further infinite families of connected Cayley graphs that have no Hamilton decomposition, one is based on $F8$ and the other on $F16$. Specifically, the families are
$K(mF8)$ and $K(mF16)$ for each positive integer $m\equiv 2\md 4$. There is a regular action of the symmetric group $\sym(4)$ on the arcs of $F8$, and so for each positive integer $m$ we have that 
$K(mF8)$ is a Cayley graph on $\sym(4)\times\Z_m$. Similarly, 
$K(mF16)$ is a Cayley graph on $\gl(2,3)\times\Z_m$, where $\gl(2,3)$ denotes the general linear group of invertible $2$ by $2$ matrices over a field with three elements. Explicitly, in the case $m=1$ we have 
$$K(F8)\cong\cay(\sym(4)\,;\{(1\ 2),(2\ 3\ 4),(2\ 4\ 3)\})$$
and 
$$K(F16)\cong\cay(\gl(2,3)\,;\{A,B,B^{-1}\}){\rm ~where~}
A=\left( \begin{array}{cc}
0 & 1 \\
1 & 0 \\
\end{array} \right)
{\rm ~and~}
B=\left( \begin{array}{cc}
1 & 1 \\
0 & 1 \\
\end{array} \right).$$

To see that 
$K(mF8)$ has no Hamilton decomposition when $m\equiv 2\md 4$,
first observe that $F8$ contains only six distinct Hamilton cycles. 
Let these Hamilton cycles be $Y_1,Y_2,\ldots,Y_6$. Also, for $1\leq i\leq 6$, let $n_i$ be the number of copies of $Y_i$ in a putative Hamilton decomposition of $mF8$.
If $u$ and $v$ are adjacent vertices in $mF8$, then it follows that the equation
$\Sigma_{i=1}^6\delta_in_i=m$ holds, where $\delta_i=1$ if $Y_i$ has an edge with endpoints $u$ and $v$, and $\delta_i=0$ otherwise. 
The twelve edges of $F8$ thus give us twelve equations in the variables
$n_1,n_2,\ldots,n_6$, and it is routine to check that these have no integral solution
when $m\equiv 2\md 4$.
It follows that $mF8$ has no Hamilton decomposition when $m\equiv 2\md 4$. 
So applying Lemma \ref{mainconstruction} gives us the following result.

\begin{theorem}\label{CayleyGraphsCube}
For each positive integer $m\equiv 2\md 4$, 
$K(mF8)$ is a connected Cayley graph 
that has no Hamilton decomposition.
\end{theorem}

Using similar arguments it can also be shown that $K(mF16)$ also has no Hamilton decomposition when 
$m\equiv 2\md 4$, which gives us the following theorem.

\begin{theorem}\label{CayleyGraphsF16}
For each positive integer $m\equiv 2\md 4$, $K(mF16)$ is a connected Cayley graph 
that has no Hamilton decomposition.
\end{theorem}

\section{Concluding remarks and questions}

In Section \ref{6valentSection} we mentioned
a computer check for the existence of Hamilton decompositions of $2X$, where $X$ is a $3$-valent arc-transitive graph of order at most $50$.
We have also verified by computer whether there exists a Hamilton decomposition of $3X$ for each $3$-valent arc-transitive graph $X$ of order at most $50$.
Every such graph has a Hamilton decomposition, except that $3F10$, $3F24$ and $3F28$ have no Hamilton decomposition.
Thus, by Lemma \ref{mainconstruction}, the graphs
$K(3F10)$, $K(3F24)$ and $K(3F28)$
are connected $9$-valent vertex-transitive graphs that have
no Hamilton decomposition. The fact that $K(3F10)$ and $K(3F28)$
have no Hamilton decomposition has been noted previously.

We now know that for infinitely many values of $k$, including $k=3$, $4$ and $6$, 
there exist connected $k$-valent vertex-transitive graphs that have no Hamilton decomposition. 
It is natural to ask whether such graphs exist for all $k\geq 3$. The smallest undecided valency is $k=5$.
One may ask the same question in relation to connected Cayley graphs. 
However, in this case the smallest undecided valency is $k=3$.  
Indeed, it is a well-known conjecture that all connected Cayley graphs 
have a Hamilton cycle (except the complete graph of order $2$), which of course 
implies that all connected $3$-valent Cayley graphs have a Hamilton decomposition. 

The graph $K(2F8)$ is a connected Cayley graph of order $48$ that
has no Hamilton decomposition.
It would be interesting to know if there exist any connected Cayley graphs of order less than $48$
that have no Hamilton decomposition. Any such graph has order at least $32$. 
It would also be interesting to know whether $K(2F8)$ is the smallest 
connected $6$-valent vertex-transitive graph that has no Hamilton decomposition, and whether 
$K(L(F10))$ is the smallest connected $4$-valent Cayley graph that has no Hamilton decomposition. 
Another open question is whether there are any connected Cayley graphs of odd order that have no Hamilton decomposition. At present, $L(F10)$ is the only connected vertex-transitive graph of odd order that is known to have no Hamilton decomposition.

\vspace{0.2cm} \noindent{\bf Acknowledgement.}
The authors acknowledge the support of the Australian Research Council via grants DP120100790 and DP120103067, and are grateful to Marston Conder for helpful discussions.


\vspace{-0.55cm}

\begin{thebibliography}{99}

{\small

\bibitem{Als1}
B. Alspach, Research Problem 59, {\it Discrete Math.}, {\bf 50} (1984), 115.

\bibitem{Als2}
B. Alspach, Unsolved Problem 4.5, {\it Ann. Discrete Math.}, 
{\bf 27} (1985), 464.

\bibitem{Als3}
B. Alspach,
Hamiltonian partitions of vertex-transitive graphs of order $2p$, 
Proceedings of the Eleventh Southeastern Conference on Combinatorics, 
Graph Theory and Computing (Florida Atlantic Univ., Boca Raton, Fla., 1980), 
Vol. I.  
{\it Congr. Numer.},
{\bf 28} (1980), 217--221. 

\bibitem{AlsBryKre}
B. Alspach, D. Bryant and D. L. Kreher, 
Vertex-Transitive Graphs Of Prime-Squared Order Are Hamilton-Decomposable, 
{\it J. Combin. Des.}, 
{\bf 22} (2014), 12--25.

\bibitem{AlsDob}
B. Alspach and E. Dobson,
On automorphism groups of graph truncations,
Ars Mathematica Contemporanea,
(to appear).

\bibitem{AlsHeiLiu}
B. Alspach, K. Heinrich and G. Z. Liu, 
Orthogonal factorizations of graphs,
Contemporary design theory, 13--40, Wiley-Intersci. Ser. Discrete Math. Optim., Wiley, New York, 1992.

\bibitem{AlsRos}
B. Alspach and M. Rosenfeld, 
On Hamilton decompositions of prisms over simple $3$-polytopes,
{\it Graphs Combin.},
{\bf 2} (1986), 1--8. 

\bibitem{Bab1}
L. Babai,
Problem 17, Unsolved Problems, in: Summer Research Workshop in Algebraic Combinatorics, 
Simon Fraser University, July 1979.

\bibitem{Bab2}
L. Babai,
Automorphism groups, Isomorphism, Reconstruction, in:
R. L. Graham, M. Grotschel, L. Lov\'asz, (Eds.),
Handbook of Combinatorics, North--Holland, 1995, pp 1447--1540 (Chapter 27).

\bibitem{Ber}
J-C. Bermond,
Hamiltonian graphs, in 
L. W. Beineke and R. J. Wilson (Eds.),
Selected Topics in Graph Theory,
Academic Press, London, 1978, pp. 127--167.

\bibitem{BerFavMah} J-C. Bermond, O. Favaron, and M. Mah\'{e}o, 
Hamiltonian decomposition of Cayley graphs of degree $4$, 
{\it J. Combin. Theory Ser. B}, {\bf 46} (1989), 142--153.

\bibitem{CheQui}
C. C. Chen and N. F. Quimpo, 
On strongly Hamiltonian abelian group graphs. Combinatorial mathematics, VIII (Geelong, 1980), pp. 23--34, Lecture Notes in Math., 884, Springer, Berlin-New York, 1981. 


\bibitem{ConDob}
M. Conder and P. Dobcs\'anyi, 
Trivalent symmetric graphs on up to 768 vertices,
{\it J. Combin. Math. Combin. Comput.},
{\bf 40} (2002), 41--63.

\bibitem{ConNed}
M. Conder and R. Nedela, 
A refined classification of symmetric cubic graphs,
{\it J. Algebra},
{\bf 322} (2009), 722--740.

\bibitem{CsaKuhLoOstTre}
B. Csaba, D. K\"uhn, A. Lo, D. Osthus and A. Treglown, 
Proof of the 1-factorization and Hamilton Decomposition Conjectures,
{\it Mem. Amer. Math. Soc.}, 
(to appear).

\bibitem{Dea1}
M. Dean,
On Hamilton cycle decomposition of $6$-regular circulant graphs,
{\it Graphs Combin.},
{\bf 22} (2006), 331--340. 

\bibitem{Dea2}
M. Dean, 
Hamilton cycle decomposition of $6$-regular circulants of odd order, 
{\it J. Combin. Des.},
{\bf 15} (2007), 91--97.

\bibitem{FanLicLiu}
C. Fan, D. R. Lick and J. Liu, 
Pseudo-Cartesian product and Hamiltonian decompositions of Cayley graphs on abelian groups,
{\it Discrete Math.},
{\bf 158} (1996), 49--62. 

\bibitem{Gou}
R. J. Gould, 
Advances on the Hamiltonian problem -- a survey,
{\it Graphs Combin.},\
{\bf 19} (2003), 7--52.

\bibitem{GruMal}
B. Gr\"unbaum and J. Malkevitch, 
Pairs of edge-disjoint Hamiltonian circuits,
{\it Aequationes Math.}
{\bf 14} (1976), 191--196. 

\bibitem{Jac}
B. Jackson,
Edge-disjoint Hamilton cycles in regular graphs of large degree,
{\it J. London Math. Soc. (2)}
{\bf 19} (1979), 13--16.

\bibitem{Kot}
A. Kotzig,
Hamilton graphs and Hamilton circuits. 
{\it Theory of Graphs and its Applications (Proc. Sympos. Smolenice 1963)}, 
Nakl. CSAV, Praha: 63--82, 1964.

\bibitem{KutMar}
K. Kutnar and D. Maru\u si\u c, 
Hamilton cycles and paths in vertex-transitive graphs -- current directions,
{\it Discrete Math.},
{\bf 309} (2009), 5491--5500.

\bibitem{Lau}
P. J. Laufer,
On strongly Hamiltonian complete bipartite graphs,
{\it Ars Combin.}, 
{\bf 9} (1980), 43--46.

\bibitem{Liu1} 
J. Liu, 
Hamiltonian decompositions of Cayley graphs on abelian groups, 
{\it Discrete Math.}, 
{\bf 131} (1994), 163--171.

\bibitem{Liu2}
J. Liu, 
Hamiltonian decompositions of Cayley graphs on abelian groups of odd order, 
{\it J. Combin. Theory Ser. B},
{\bf 66} (1996), 75--86.

\bibitem{Liu3}
J. Liu, 
Hamiltonian decompositions of Cayley graphs on abelian groups of even order, 
{\it J. Combin. Theory Ser. B},
{\bf 88} (2003), 305--321.

\bibitem{Lov}
L. Lov\'asz, 
Combinatorial structures and their applications, in: Proc. Calgary Internat. Conf. Calgary, Alberta, 1969, Gordon and Breach, New York, 1970, pp. 243-246, Problem 11.

\bibitem{Mad}
W. Mader, Minimale $n$-fach kantenzusammenh\"angende Graphen, 
{\it Math. Ann.}, 
{\bf 191} (1971), 21--28.

\bibitem{McKRoy}
B. D. McKay and G. F. Royle,
The transitive graphs with at most $26$ vertices,
{\it Ars Combin.},
{\bf 30} (1990), 161--176.

\bibitem{Nas}
C. St. J. A. Nash-Williams,
Hamiltonian arcs and circuits. 1971 Recent Trends in Graph Theory (Proc. Conf., New York, 1970) pp. 197--210 Lecture Notes in Mathematics, Vol. 186 Springer, Berlin.

\bibitem{Pik1}
D. A. Pike, 
Hamilton decompositions of line graphs of perfectly 1-factorisable graphs of even degree, {\it Australas. J. Combin.},
{\bf 12} (1995), 291--294.

\bibitem{Pik2}
D. A. Pike, 
Snarks and non-Hamiltonian cubic 2-edge-connected graphs of small order,
{\it J. Combin. Math. Combin. Comput.},
{\bf 23} (1997), 129--141.

\bibitem{PotSpiVer}
P. Poto\u cnik, P. Spiga and G. Verret, 
Cubic vertex-transitive graphs on up to 1280 vertices, 
{\it J. Symbolic Comput.},
{\bf 50} (2013), 465--477.

\bibitem{Wei}
E. W. Weisstein, 
Hamilton Decomposition,
From {\it MathWorld} -- A Wolfram Web Resource. http://mathworld.wolfram.com/HamiltonDecomposition.html

\bibitem{Wes1} 
E. E. Westlund, 
Hamilton decompositions of $6$-regular Cayley graphs on even Abelian groups with involution-free connections sets,
{\it Discrete Math.},
{\bf 331} (2014), 117--132.

\bibitem{Wes2} 
E. E. Westlund, 
Hamilton decompositions of certain $6$-regular Cayley graphs on Abelian groups with a cyclic subgroup of index two, 
{\it Discrete Math.},
{\bf 312} (2012), 3228--3235. 

\bibitem{WesLiuKre}
E. E. Westlund, J. Liu and D. L. Kreher, 
$6$-regular Cayley graphs on abelian groups of odd order are Hamiltonian decomposable,
{\it Discrete Math.},
{\bf 309} (2009), 5106--5110. 


}
\end{thebibliography}
\end{document}